\documentclass[12pt]{article}

\usepackage{amsmath,amssymb,amsthm,amscd,a4wide,upref}

\usepackage{hyperref}
\hypersetup{colorlinks,citecolor=blue,filecolor=black,linkcolor=blue,urlcolor=blue}
\usepackage{enumerate}
\usepackage{mathtools}

\usepackage{lmodern}     % set math font to Latin modern math
\usepackage[T1]{fontenc}
\begin{document}

% % % % % % % % % % % % % %
%\newcommand{\tma}{\textcolor{magenta}}
%\newcommand{\tre}{\textcolor{red}}
% % % % % % % % % % % % % %

\newcommand{\ad}{{\rm ad}}
\newcommand{\cri}{{\rm cri}}
\newcommand{\row}{{\rm row}}
\newcommand{\col}{{\rm col}}
\newcommand{\End}{{\rm{End}\ts}}
\newcommand{\Rep}{{\rm{Rep}\ts}}
\newcommand{\Hom}{{\rm{Hom}}}
\newcommand{\Det}{{\rm{Det}}}
\newcommand{\Mat}{{\rm{Mat}}}
\newcommand{\ch}{{\rm{ch}\ts}}
\newcommand{\chara}{{\rm{char}\ts}}
\newcommand{\diag}{{\rm diag}}
\newcommand{\non}{\nonumber}
\newcommand{\wt}{\widetilde}
\newcommand{\wh}{\widehat}
\newcommand{\ol}{\overline}
\newcommand{\ot}{\otimes}
\newcommand{\la}{\lambda}
\newcommand{\La}{\Lambda}
\newcommand{\De}{\Delta}
\newcommand{\al}{\alpha}
\newcommand{\be}{\beta}
\newcommand{\ga}{\gamma}
\newcommand{\Ga}{\Gamma}
\newcommand{\ep}{\epsilon}
\newcommand{\ka}{\kappa}
\newcommand{\vk}{\varkappa}
\newcommand{\si}{\sigma}
\newcommand{\vs}{\varsigma}
\newcommand{\vp}{\varphi}
\newcommand{\vpi}{\varpi}
\newcommand{\de}{\delta}
\newcommand{\ze}{\zeta}
\newcommand{\om}{\omega}
\newcommand{\Om}{\Omega}
\newcommand{\ee}{\epsilon^{}}
\newcommand{\su}{s^{}}
\newcommand{\hra}{\hookrightarrow}
\newcommand{\ve}{\varepsilon}
\newcommand{\ts}{\,}
\newcommand{\vac}{\mathbf{1}}
\newcommand{\di}{\partial}
\newcommand{\qin}{q^{-1}}
\newcommand{\tss}{\hspace{1pt}}
\newcommand{\pr}{^{\tss\prime}}
\newcommand{\tra}{ {\rm t}}
\newcommand{\Sr}{ {\rm S}}
\newcommand{\U}{ {\rm U}}
\newcommand{\BL}{ {\overline L}}
\newcommand{\BE}{ {\overline E}}
\newcommand{\BP}{ {\overline P}}
\newcommand{\AAb}{\mathbb{A}\tss}
\newcommand{\CC}{\mathbb{C}\tss}
\newcommand{\KK}{\mathbb{K}\tss}
\newcommand{\QQ}{\mathbb{Q}\tss}
\newcommand{\SSb}{\mathbb{S}\tss}
\newcommand{\TT}{\mathbb{T}\tss}
\newcommand{\ZZ}{\mathbb{Z}\tss}
\newcommand{\DY}{ {\rm DY}}
\newcommand{\X}{ {\rm X}}
\newcommand{\Y}{ {\rm Y}}
\newcommand{\Z}{{\rm Z}}
\newcommand{\Ac}{\mathcal{A}}
\newcommand{\Dc}{\mathcal{D}}
\newcommand{\Lc}{\mathcal{L}}
\newcommand{\Mc}{\mathcal{M}}
\newcommand{\Pc}{\mathcal{P}}
\newcommand{\Qc}{\mathcal{Q}}
\newcommand{\Rc}{\mathcal{R}}
\newcommand{\Sc}{\mathcal{S}}
\newcommand{\Tc}{\mathcal{T}}
\newcommand{\Bc}{\mathcal{B}}
\newcommand{\Ec}{\mathcal{E}}
\newcommand{\Fc}{\mathcal{F}}
\newcommand{\Gc}{\mathcal{G}}
\newcommand{\Hc}{\mathcal{H}}
\newcommand{\Uc}{\mathcal{U}}
\newcommand{\Vc}{\mathcal{V}}
\newcommand{\Wc}{\mathcal{W}}
\newcommand{\Yc}{\mathcal{Y}}
\newcommand{\Ar}{{\rm A}}
\newcommand{\Br}{{\rm B}}
\newcommand{\Ir}{{\rm I}}
\newcommand{\Fr}{{\rm F}}
\newcommand{\Jr}{{\rm J}}
\newcommand{\Or}{{\rm O}}
\newcommand{\GL}{{\rm GL}}
\newcommand{\Spr}{{\rm Sp}}
\newcommand{\Rr}{{\rm R}}
\newcommand{\Zr}{{\rm Z}}
\newcommand{\gl}{\mathfrak{gl}}
\newcommand{\middd}{{\rm mid}}
\newcommand{\ev}{{\rm ev}}
\newcommand{\Pf}{{\rm Pf}}
\newcommand{\Norm}{{\rm Norm\tss}}
\newcommand{\oa}{\mathfrak{o}}
\newcommand{\spa}{\mathfrak{sp}}
\newcommand{\osp}{\mathfrak{osp}}
\newcommand{\f}{\mathfrak{f}}
\newcommand{\g}{\mathfrak{g}}
\newcommand{\h}{\mathfrak h}
\newcommand{\n}{\mathfrak n}
\newcommand{\z}{\mathfrak{z}}
\newcommand{\Zgot}{\mathfrak{Z}}
\newcommand{\p}{\mathfrak{p}}
\newcommand{\sll}{\mathfrak{sl}}
\newcommand{\agot}{\mathfrak{a}}
\newcommand{\qdet}{ {\rm qdet}\ts}
\newcommand{\Ber}{ {\rm Ber}\ts}
\newcommand{\HC}{ {\mathcal HC}}
\newcommand{\cdet}{{\rm cdet}}
\newcommand{\rdet}{{\rm rdet}}
\newcommand{\tr}{ {\rm tr}}
\newcommand{\gr}{ {\rm gr}\ts}
\newcommand{\str}{ {\rm str}}
\newcommand{\loc}{{\rm loc}}
\newcommand{\Gr}{{\rm G}}
\newcommand{\sgn}{ {\rm sgn}\ts}
\newcommand{\sign}{{\rm sgn}}
\newcommand{\ba}{\bar{a}}
\newcommand{\bb}{\bar{b}}
\newcommand{\bi}{\bar{\imath}}
\newcommand{\bj}{\bar{\jmath}}
\newcommand{\bk}{\bar{k}}
\newcommand{\bl}{\bar{l}}
\newcommand{\hb}{\mathbf{h}}
\newcommand{\Sym}{\mathfrak S}
\newcommand{\fand}{\quad\text{and}\quad}
\newcommand{\Fand}{\qquad\text{and}\qquad}
\newcommand{\For}{\qquad\text{or}\qquad}
\newcommand{\for}{\quad\text{or}\quad}
\newcommand{\grpr}{{\rm gr}^{\tss\prime}\ts}
\newcommand{\degpr}{{\rm deg}^{\tss\prime}\tss}
\newcommand{\stir}{\genfrac{\{}{\}}{0pt}{}}
\newcommand{\stirone}{\genfrac{[}{]}{0pt}{}}

\renewcommand{\theequation}{\arabic{section}.\arabic{equation}}

\numberwithin{equation}{section}

\newtheorem{thm}{Theorem}[section]
\newtheorem{lem}[thm]{Lemma}
\newtheorem{prop}[thm]{Proposition}
\newtheorem{cor}[thm]{Corollary}
\newtheorem{conj}[thm]{Conjecture}
\newtheorem*{mthm}{Main Theorem}
\newtheorem*{mthma}{Theorem A}
\newtheorem*{mthmb}{Theorem B}
\newtheorem*{mthmc}{Theorem C}
\newtheorem*{mthmd}{Theorem D}

\theoremstyle{definition}
\newtheorem{defin}[thm]{Definition}

\theoremstyle{remark}
\newtheorem{remark}[thm]{Remark}
\newtheorem{example}[thm]{Example}
\newtheorem{examples}[thm]{Examples}

\newcommand{\bth}{\begin{thm}}
\renewcommand{\eth}{\end{thm}}
\newcommand{\bpr}{\begin{prop}}
\newcommand{\epr}{\end{prop}}
\newcommand{\ble}{\begin{lem}}
\newcommand{\ele}{\end{lem}}
\newcommand{\bco}{\begin{cor}}
\newcommand{\eco}{\end{cor}}
\newcommand{\bde}{\begin{defin}}
\newcommand{\ede}{\end{defin}}
\newcommand{\bex}{\begin{example}}
\newcommand{\eex}{\end{example}}
\newcommand{\bes}{\begin{examples}}
\newcommand{\ees}{\end{examples}}
\newcommand{\bre}{\begin{remark}}
\newcommand{\ere}{\end{remark}}
\newcommand{\bcj}{\begin{conj}}
\newcommand{\ecj}{\end{conj}}

\newcommand{\bal}{\begin{aligned}}
\newcommand{\eal}{\end{aligned}}
\newcommand{\beq}{\begin{equation}}
\newcommand{\eeq}{\end{equation}}
\newcommand{\ben}{\begin{equation*}}
\newcommand{\een}{\end{equation*}}

\newcommand{\bpf}{\begin{proof}}
\newcommand{\epf}{\end{proof}}

\def\beql#1{\begin{equation}\label{#1}}

%%% Slaven's notation

%\newcommand{\End}{\mathop{\mathrm{End}}}
%\newcommand{\Hom}{\mathop{\mathrm{Hom}}}
%\newcommand{\vac}{\mathop{\mathrm{\boldsymbol{1}}}}
\newcommand{\Res}{\mathop{\mathrm{Res}}}

\title{\Large\bf On Segal--Sugawara vectors and Casimir elements\\
 for classical Lie algebras}

\author{A. I. Molev}

\date{} % Start August 2020
\maketitle

%\vspace{4 mm}

\begin{abstract}
We consider the centers of
the affine vertex algebras at the critical level
associated with simple Lie algebras.
We derive new formulas for generators of the centers in
the classical types.
We also give a new formula for the Capelli-type determinant
for the symplectic Lie algebras and
calculate the Harish-Chandra images
of the Casimir elements
arising from the characteristic polynomial
of the matrix of generators of each classical Lie algebra.

\end{abstract}

%\vspace{5 mm}
%%%
%%%{\it Key words:}
%%%

%\newpage

%\tableofcontents
%
%\newpage

\section{Introduction}
\label{sec:int}

Let $\g$ be a simple Lie algebra over $\CC$ equipped with
a standard symmetric invariant bilinear form.
The affine Kac--Moody algebra $\wh\g$
is defined as the central
extension
\beql{km}
\wh\g=\g\tss[t,t^{-1}]\oplus\CC K
\eeq
of the Lie algebra of Laurent
polynomials in $t$. The vacuum module $V_{\cri}(\g)$
at the critical level over $\wh\g$ is
the quotient of the
universal enveloping algebra $\U(\wh\g)$ by the left
ideal generated by $\g[t]$ and $K+h^{\vee}$, where $h^{\vee}$ is the dual Coxeter number for $\g$.
The vacuum module has a vertex algebra structure and is known as
the {\em {\em({\em universal})} affine vertex algebra}; see e.g.
\cite{f:lc} and \cite{k:va} for definitions.
The {\em center} of the vertex algebra $V_{\cri}(\g)$
is defined by
\ben
\z(\wh\g)
=\{S\in V_{\cri}(\g)\ |\ \g[t]\ts S=0\}.
\een
Any element of $\z(\wh\g)$
is called a {\it Segal--Sugawara vector\/}.
The vertex algebra axioms imply that the center is a commutative
associative algebra which can be regarded as a subalgebra of $\U\big(t^{-1}\g[t^{-1}]\big)$.
The algebra $\z(\wh\g)$ is equipped with the derivation $T=-d/dt$
arising from the vertex algebra structure.
By a theorem of Feigin and Frenkel~\cite{ff:ak}, the differential algebra
$\z(\wh\g)$ possesses generators
$S_1,\dots,S_n$ so that
$\z(\wh\g)$ is the algebra of polynomials
\ben
\z(\wh\g)=\CC[T^{\tss r}S_l\ |\ l=1,\dots,n,\ \ r\geqslant 0],
\een
where $n=\text{rank}\ts\g$; see also \cite{f:lc}. The algebra $\z(\wh\g)$
is known as the {\em Feigin--Frenkel center}, and we call $S_1,\dots,S_n$
a {\em complete set of Segal--Sugawara vectors}. According to \cite{ff:ak}
(see also \cite{f:lc}),
the center can be identified
with the {\em classical $\Wc$-algebra}
associated with the Langlands dual Lie algebra ${}^L\g$
via an affine version of the Harish-Chandra isomorphism
\beql{hchiaff}
\z(\wh\g)\cong \Wc({}^L\g).
\eeq

Explicit formulas for complete sets of Segal--Sugawara vectors
were given in \cite{cm:ho} and \cite{ct:qs} for the Lie algebras $\g$ of type $A$,
and in \cite{m:ff} for types $B$, $C$ and $D$ with the use of the Brauer algebra.
The images of the vectors with respect to the Harish-Chandra
isomorphism \eqref{hchiaff} are easy to find in type $A$ \cite{cm:ho},
whereas the calculations in the remaining classical types
in \cite{mm:yc}
rely on the $q$-character formulas for certain Yangian representations.
A detailed exposition of these results together with applications to
commutative subalgebras in enveloping algebras and to higher order Hamiltonians
in the Gaudin models can be found in \cite{m:so}.
A complete set of Segal--Sugawara vectors
for the Lie algebra of type $G_2$ was produced in \cite{mrr:ss}
by using computer-assisted calculations.
A different method to construct generators of $\z(\wh\g)$ was developed
in \cite{o:sf} which lead to new explicit formulas in the case of the Lie algebras
of types $B,C,D$ and $G_2$.

In this paper we derive new uniform expressions for the Segal--Sugawara vectors
in all classical types. In types $B$, $C$ and $D$
we transform the formulas produced in \cite{m:ff} by eliminating
the dependence on the Brauer diagrams with horizontal edges. In particular, the vectors are
given more explicitly in the symplectic case thus resolving the `analytic continuation'
procedure used in \cite{m:ff}; see also \cite[Ch.~8]{m:so}.
We also show that both in the orthogonal and symplectic case
the Segal--Sugawara vectors of \cite{m:ff} coincide with those
in \cite{o:sf}.

In all classical types we also consider the Casimir elements
obtained by the application of
the symmetrization map to basic $\g$-invariants in the symmetric algebra $\Sr(\g)$,
arising from the characteristic polynomial
of the matrix of generators. We
calculate their Harish-Chandra images
in terms of shifted invariant polynomials.

Certain analogues of the Capelli determinant for the
orthogonal and symplectic Lie algebras were obtained in \cite{m:sd}
with the use of the Sklyanin determinants; see also \cite{hu:ci} and \cite{w:ce}
for different approaches. Here we give a new symmetrized determinant
expression for the Capelli-type determinant in the symplectic case.

\section{Segal--Sugawara vectors}
\label{sec:ssv}

In all classical types, the new formulas for Segal--Sugawara vectors will take the form of
linear combinations of certain symmetrized $\la$-minors
or $\la$-permanents associated with partitions $\la$.
We consider the Lie algebras of type $A$ first, where
the formulas are derived from the results of
\cite{cm:ho} and \cite{ct:qs}; see also \cite[Ch.~7]{m:so}.

\subsection{Generators of $\z(\wh\gl_N)$}
\label{subsec:ga}

The affine Kac--Moody algebra $\wh\gl_N=\gl_N[t,t^{-1}]\oplus\CC K$
has the commutation relations
\beql{commrel}
\big[E_{ij}[r],E_{kl}[s\tss]\tss\big]
=\de_{kj}\ts E_{i\tss l}[r+s\tss]
-\de_{i\tss l}\ts E_{kj}[r+s\tss]
+r\tss\de_{r,-s}\ts K\Big(\de_{kj}\tss\de_{i\tss l}
-\frac{\de_{ij}\tss\de_{kl}}{N}\Big),
\eeq
and the element $K$ is central. Here and below we write $X[r]$ for the
element $Xt^r$ of the Lie algebra of Laurent polynomials $\g[t,t^{-1}]$
with $X\in\g$ and $r\in\ZZ$. The elements $E_{ij}$ form a standard basis of $\gl_N$.
The
critical level $K=-N$ coincides with
the negative of the dual Coxeter number for $\sll_N$.

Let $\la=(\la_1,\dots,\la_{\ell})$ be a partition of $m$ of length $\ell=\ell(\la)$,
so that $\la_1\geqslant\dots\geqslant\la_{\ell}>0$ and $\la_1+\dots+\la_{\ell}=m$.
We will denote by $c^{}_{\la}$ the number
of permutations in the symmetric group $\Sym_m$ of cycle type $\la$; see \eqref{cla} below
for an explicit expression.
The
{\em symmetrized $\la$-minors} $D(\la)$ and {\em symmetrized $\la$-permanents} $P(\la)$
are elements of $V_{\cri}(\gl_N)\cong\U\big(t^{-1}\gl_N[t^{-1}]\big)$ defined by
\begin{align}
\non
D(\la)&=
\frac{1}{{\ell\tss}!}\ts \sum_{i_1,\dots,\ts i_{\ell}=1}^N\ts
\sum_{\si\in\Sym_{\ell}}\sgn\si\cdot
E_{i_{\si(1)}\tss i_{1}}[-\la_1]\dots E_{i_{\si({\ell})}\tss i_{{\ell}}}[-\la_{\ell}]\\
\intertext{and}
P(\la)&=
\frac{1}{{\ell\tss}!}\ts \sum_{i_1,\dots,\ts i_{\ell}=1}^N\ts
\sum_{\si\in\Sym_{\ell}}
E_{i_{\si(1)}\tss i_{1}}[-\la_1]\dots E_{i_{\si({\ell})}\tss i_{{\ell}}}[-\la_{\ell}].
\non
\end{align}

\bth\label{thm:typea}
All elements
\ben
\phi^{}_{m}=\sum_{\la\ts\vdash m}\ts \binom{N}{{\ell}}^{-1}\ts
c^{}_{\la}\ts D(\la)
\Fand
\psi^{}_{m}=\sum_{\la\ts\vdash m}\ts \binom{N+\ell-1}{{\ell}}^{-1}\ts
c^{}_{\la}\ts P(\la)
\een
belong to the Feigin--Frenkel center $\z(\wh\gl_N)$. Moreover,
each family $\phi^{}_1,\dots,\phi^{}_N$ and $\psi^{}_1,\dots,\psi^{}_N$ is
a complete set of Segal--Sugawara vectors for $\gl_N$.
\eth

\bpf
The theorem will follow from the relations
\beql{phipsi}
\phi^{}_{m\tss m}=\binom{N}{m}\ts\phi^{}_{m}\Fand
\psi^{}_{m\tss m}=\binom{N+m-1}{m}\ts\psi^{}_{m}
\eeq
for the Segal--Sugawara vectors $\phi^{}_{m\tss m}$ and $\psi^{}_{m\tss m}$
used in \cite[Ch.~7]{m:so}.
To make the connection,
for any $r\in\ZZ$ combine the elements $E_{ij}[r]$
into the matrix $E[r]$ so that
\ben
E[r]=\sum_{i,j=1}^N e_{ij}\ot E_{ij}[r]\in \End\CC^N\ot \U(\wh\gl_N),
\een
where the $e_{ij}$ denote the standard matrix units.
For
each $a\in\{1,\dots,m\}$
introduce the element $E[r]_a$ of the algebra
\beql{tenprka}
\underbrace{\End\CC^{N}\ot\dots\ot\End\CC^{N}}_m{}\ot\U(\wh\gl_N)
\eeq
by
\ben
E[r]_a=\sum_{i,j=1}^{N}
1^{\ot(a-1)}\ot e_{ij}\ot 1^{\ot(m-a)}\ot E_{ij}[r].
\een
The symmetric group $\Sym_m$ acts on the space
\beql{tenprk}
(\CC^N)^{\ot m}=\underbrace{\CC^{N}\ot\CC^{N}\ot\dots\ot\CC^{N}}_m
\eeq
by permuting the tensor factors.
Denote by $H^{(m)}$ and $A^{(m)}$ the elements of the algebra
\eqref{tenprka} (with the identity components in $\U(\wh\gl_N)$) which are
the respective images of the symmetrizer $h^{(m)}$ and anti-symmetrizer $a^{(m)}$
defined by
\beql{ha}
h^{(m)}=\frac{1}{m!}\sum_{s\in\Sym_m} s
\Fand
a^{(m)}=\frac{1}{m!}\sum_{s\in\Sym_m} \sgn s\cdot s,
\eeq
under the action of $\Sym_m$.
By \cite[Thms~7.1.3 \& 7.1.4]{m:so}, the corresponding
claims of Theorem~\ref{thm:typea} hold for the elements
\beql{deftra}
\phi^{}_{m\tss m}=\tr^{}_{1,\dots,m}\ts A^{(m)} \big(T+E[-1]_1\big)\dots \big(T+E[-1]_m\big)\tss 1
\eeq
and
\beql{deftrh}
\psi^{}_{m\tss m}=\tr^{}_{1,\dots,m}\ts H^{(m)} \big(T+E[-1]_1\big)\dots \big(T+E[-1]_m\big)\tss 1,
\eeq
where the vacuum vector of $V_{\cri}(\gl_N)$ is identified with the element
$1\in \U\big(t^{-1}\gl_N[t^{-1}]\big)$
which is annihilated by the
derivation $T=-d/dt$. The trace is taken with respect to all $m$ copies of
the endomorphism algebra $\End\CC^N$
in \eqref{tenprka}.
Expand the product in \eqref{deftra} by using the relations $\big[T,E[-r]_a\big]=r\tss E[-r-1]_a$
so that $\phi^{}_{m\tss m}$ will take the form of a linear combination of the traces
\beql{trasfa}
\tr^{}_{1,\dots,m}\ts A^{(m)} E[-r_1]_{a_1}\dots E[-r_s]_{a_s}
\eeq
with $a_1<\dots<a_s$ and $r_i\geqslant 1$. The defining relations \eqref{commrel} imply
that for $a<b$ we have
\ben
E[-r]_a\ts E[-s]_b-E[-s]_b\ts E[-r]_a=E[-r-s]_a\tss P_{a\tss b}-P_{a\tss b}\ts E[-r-s]_a,
\een
where $P_{a\tss b}$ is the permutation operator
\beql{pdef}
P_{a\tss b}=\sum_{i,j=1}^N 1^{\ot(a-1)}\ot e_{ij}
\ot 1^{\ot(b-a-1)}\ot e_{ji}\ot 1^{\ot(m-b)}.
\eeq
Hence, by the cyclic property of trace,
any permutation of the factors $E[-r_i]_{a_i}$
in the expression
\eqref{trasfa}
does not change its value.
Therefore,
applying conjugations by
suitable permutations of the index set $1,\dots,m$ and using the cyclic property
of trace, we can write
\beql{expa}
\tr^{}_{1,\dots,m}\ts A^{(m)} \big(T+E[-1]_1\big)\dots \big(T+E[-1]_m\big)\tss 1
=\tr^{}_{1,\dots,m}\ts A^{(m)} \sum_{\la\ts\vdash m}\tss
c^{}_{\la}\ts E[-\la],
\eeq
for certain nonnegative integers $c^{}_{\la}$,
where we set
\ben
E[-\la]=E[-\la_1]_1 \dots E[-\la_{\ell}]_{\ell}
\een
for any partition $\la=(\la_1,\dots,\la_{\ell})$ with $\ell=\ell(\la)$.
Identifying partitions with their Young diagrams, we can write
the expression on the left hand side of \eqref{expa} as
\ben
\tr^{}_{1,\dots,m}\ts A^{(m)}\ts T\sum_{\mu\ts\vdash m-1}\tss
c^{}_{\mu}\ts E[-\mu]+
\tr^{}_{1,\dots,m}\ts A^{(m)} \sum_{\mu\ts\vdash m-1}\tss
c^{}_{\mu}\ts E[-\mu^+],
\een
where $\mu^+$ is the diagram obtained from $\mu$ by adding one box
to the first column.
Hence the coefficients $c^{}_{\la}$ satisfy the recurrence relation
\beql{rerel}
c^{}_{\la}=\sum_{\mu}\ga(\mu,\la)\tss c^{}_{\mu},
\eeq
summed over the diagrams $\mu$ which are obtained from $\la$ by removing one box, where
\ben
\ga(\mu,\la)=\begin{cases}\mu_i\cdot\text{mult}\ts(\mu_i)\qquad&\text{if}\quad\mu_i\geqslant 1,\\
                             1\qquad&\text{if}\quad\mu_i=0,
                             \end{cases}
\een
assuming the box is removed in row $i$, and $\text{mult}\ts(\mu_i)$ denotes the multiplicity
of $\mu_i$ as a part of $\mu$. Writing the partition $\la$
in the multiplicity form
$\la=(1^{\al_1}2^{\al_2}\dots m^{\al_m})$, we derive from \eqref{rerel} by induction on $m$
that
\beql{cla}
c^{}_{\la}=\frac{m!}{1^{\al_1}\al_1!\ts 2_{}^{\al_2}\al_2!\dots m^{\al_m}\al_m!}
\eeq
which is the number
of permutations in $\Sym_m$ of cycle type $\la$.

Note that the partial traces of the anti-symmetrizer are found by
\beql{patr}
\tr^{}_{\ell+1,\dots,m}\ts A^{(m)}=\binom{N}{m}\binom{N}{\ell}^{-1} A^{(\ell)},
\eeq
and so \eqref{expa} implies
\ben
\phi^{}_{m\tss m}=\binom{N}{m}\sum_{\la\ts\vdash m}\tss \binom{N}{\ell}^{-1}
c^{}_{\la}\ts \tr^{}_{1,\dots,\ell}\ts A^{(\ell)} E[-\la],
\een
which proves the first relation in \eqref{phipsi} because
\beql{dla}
D(\la)=\tr^{}_{1,\dots,\ell}\ts A^{(\ell)} E[-\la].
\eeq

The second relation in \eqref{phipsi} is verified by the same argument,
where the partial traces of the symmetrizer are evaluated by
\beql{patrh}
\tr^{}_{\ell+1,\dots,m}\ts H^{(m)}=\binom{N+m-1}{m}\binom{N+\ell-1}{\ell}^{-1} H^{(\ell)},
\eeq
and the identity
\ben
P(\la)=\tr^{}_{1,\dots,\ell}\ts H^{(\ell)} E[-\la],
\een
is used in place of \eqref{dla}.
\epf

As was pointed out in \cite[Sec.~2]{o:sf}, since the Feigin--Frenkel center $\z(\wh\gl_N)$
is invariant under the automorphism $\theta$ of $\U\big(t^{-1}\gl_N[t^{-1}]\big)$
taking $E_{ij}[r]$ to $-E_{ji}[r]$, the Segal--Sugawara vectors
of Theorem~\ref{thm:typea} can be modified to become eigenvectors of $\theta$.
Note that both $A^{(\ell)}$ and $H^{(\ell)}$ are stable under the simultaneous transpositions
with respect to all $\ell$ copies of $\End\CC^N$ and so
\ben
\theta:D(\la)\mapsto (-1)^{\ell}\tss D(\la)\Fand \theta:P(\la)\mapsto (-1)^{\ell}\tss P(\la).
\een
The same symmetry properties hold for the
anti-automorphism of $\U\big(t^{-1}\gl_N[t^{-1}]\big)$
taking $E_{ij}[r]$ to $-E_{ij}[r]$.
As in \cite[Ch.~7]{m:so}, this leads to the following.

\bco\label{cor:even}
All elements
\ben
\phi^{\circ}_{m}=\sum_{\la\ts\vdash m,\ m-\ell\ts \text{even}}\ts \binom{N}{{\ell}}^{-1}\ts
c^{}_{\la}\ts D(\la)
\Fand
\psi^{\circ}_{m}=\sum_{\la\ts\vdash m,\ m-\ell\ts \text{even}}\ts \binom{N+\ell-1}{{\ell}}^{-1}\ts
c^{}_{\la}\ts P(\la)
\een
belong to the Feigin--Frenkel center $\z(\wh\gl_N)$. Moreover,
each family $\phi^{\circ}_1,\dots,\phi^{\circ}_N$ and $\psi^{\circ}_1,\dots,\psi^{\circ}_N$ is
a complete set of Segal--Sugawara vectors for $\gl_N$.
\qed
\eco

\subsection{Generators of $\z(\wh\g)$ for types $B$, $C$ and $D$}
\label{subsec:nf}

We will regard the orthogonal Lie algebras $\oa_N$ with $N=2n+1$ and $N=2n$
and symplectic Lie algebra $\spa_N$ with $N=2n$
as subalgebras of $\gl_N$ spanned by the elements $F_{i\tss j}$,
\ben
F_{i\tss j}=E_{i\tss j}-E_{j\pr i\pr}\Fand F_{i\tss j}
=E_{i\tss j}-\ve_i\ts\ve_j\ts E_{j\pr i\pr},
\een
respectively, for $\oa_N$ and $\spa_N$, where
$i\pr=N-i+1$.
In the symplectic case we set
$\ve_i=1$ for $i=1,\dots,n$ and
$\ve_i=-1$ for $i=n+1,\dots,2n$.
As before, we will write $F_{ij}[r]=F_{ij}t^r$ with $r\in\ZZ$
for elements of the
Kac--Moody algebra $\wh\g$ for $\g=\oa_N$ or $\spa_N$,
as defined in \eqref{km}.

Let $\la=(\la_1,\dots,\la_{\ell})$ be a partition of $m$ of length $\ell=\ell(\la)$.
In the case $\g=\spa_{2n}$ we
introduce the
corresponding {\em symmetrized $\la$-minor} by
\ben
D(\la)=
\frac{1}{{\ell\tss}!}\ts \sum_{i_1,\dots,\ts i_{\ell}=1}^{2n}\ts
\sum_{\si\in\Sym_{\ell}}\sgn\si\cdot
F_{i_{\si(1)}\tss i_{1}}[-\la_1]\dots F_{i_{\si({\ell})}\tss i_{{\ell}}}[-\la_{\ell}].
\een
In the case $\g=\oa_N$
the {\em symmetrized $\la$-permanent} is defined by
\ben
P(\la)=
\frac{1}{{\ell\tss}!}\ts \sum_{i_1,\dots,\ts i_{\ell}=1}^N\ts
\sum_{\si\in\Sym_{\ell}}
F_{i_{\si(1)}\tss i_{1}}[-\la_1]\dots F_{i_{\si({\ell})}\tss i_{{\ell}}}[-\la_{\ell}].
\een
Both $D(\la)$ and $P(\la)$ are zero unless $\ell(\la)$ is even; see also \eqref{dpla} below.

We will keep using the numbers $c^{}_{\la}$ given by \eqref{cla} which count
the permutations in $\Sym_m$ of cycle type $\la$. Recall a distinguished
Segal--Sugawara vector $\Pf\tss F[-1]$ for $\g=\oa_{2n}$, which is the
(noncommutative) {\em Pfaffian}
of the matrix $F[-1]=\big[F_{ij}[-1]\big]$; see \cite{m:ff}, \cite[Sec.~8.1]{m:so}.

\bth\label{thm:typebcd}
All elements
\ben
\phi^{}_{m}=\sum_{\la\ts\vdash m,\  \ell(\la)\ts \text{even}}\ts \binom{2n+1}{{\ell}}^{-1}\ts
c^{}_{\la}\ts D(\la)
\fand
\phi^{}_{m}=\sum_{\la\ts\vdash m,\  \ell(\la)\ts \text{even}}\ts \binom{N+\ell-2}{{\ell}}^{-1}\ts
c^{}_{\la}\ts P(\la)
\een
belong to the Feigin--Frenkel center $\z(\wh\g)$ in the symplectic and orthogonal case,
respectively.
Moreover, the family
$\phi^{}_2,\phi^{}_4,\dots,\phi^{}_{2n}$
is a complete set of Segal--Sugawara vectors for
$\g=\spa_{2n}$ and $\g=\oa_{2n+1}$, while
$\phi^{}_2,\phi^{}_4,\dots,\phi^{}_{2n-2},\Pf\tss F[-1]$ is
complete set of Segal--Sugawara vectors for $\g=\oa_{2n}$.
\eth

\bpf
We will derive the theorem from the relations
\beql{phipsibcd}
\phi^{}_{m\tss m}=\binom{2n+1}{m}\ts\phi^{}_{m}\Fand
\phi^{}_{m\tss m}=\binom{N+m-2}{m}\ts\phi^{}_{m}
\eeq
in the symplectic and orthogonal case, respectively, for
the Segal--Sugawara vectors $\phi^{}_{m\tss m}$ considered
in \cite[Ch.~8]{m:so}.
We will regard the $N\times N$ matrix $F[r]=\big[F_{ij}[r]\big]$
as the element
\ben
F[r]=\sum_{i,j=1}^N e_{ij}\ot F_{ij}[r]\in \End\CC^N\ot \U(\wh\g).
\een
It has the skew-symmetry property $F[r]+F[r]^{\tra}=0$
with respect to the transposition defined by
\beql{transpmu}
\tra:e_{ij}\mapsto\begin{cases} e_{j'i'}\qquad&\text{in the orthogonal case,}\\
\ve_i\tss\ve_j\ts e_{j'i'}\qquad&\text{in the symplectic case.}
\end{cases}
\eeq
For
each $a\in\{1,\dots,m\}$
introduce the element $F[r]_a$ of the algebra
\beql{tenprkao}
\underbrace{\End\CC^{N}\ot\dots\ot\End\CC^{N}}_m{}\ot\U(\wh\g)
\eeq
by
\ben
F[r]_a=\sum_{i,j=1}^{N}
1^{\ot(a-1)}\ot e_{ij}\ot 1^{\ot(m-a)}\ot F_{ij}[r].
\een
The $a$-th partial transposition $\tra_a$ on the algebra \eqref{tenprkao} acts as
the map \eqref{transpmu} on the $a$-th copy of $\End\CC^N$ and as the identity map
on all other tensor factors.

The Segal--Sugawara vectors $\phi^{}_{m\tss m}$ are constructed with the use of the Brauer
algebra $\Bc_m(\om)$ whose definition we will now recall.
Consider an $m$-{\em diagram}
$d$ which is a collection
of $2m$ dots arranged into two rows with $m$ dots in each row labelled by $1,\dots,m$;
the dots are
connected by $m$ edges in such a way that any dot belongs to only one edge.
An edge is called {\em horizontal} if it connects two dots in the same row.
The product $d\tss d^{\tss\prime}$ of two diagrams $d$ and $d^{\tss\prime}$ is determined by
placing $d$ under $d^{\tss\prime}$ and identifying the vertices
of the bottom row of $d^{\tss\prime}$ with the corresponding
vertices in the top row of $d$. Let $s$ be the number of
closed loops obtained in this placement. The product $d\tss d^{\tss\prime}$ is given by
$\om^{\tss s}$ times the resulting diagram without loops.
The algebra
$\Bc_m(\om)$ is defined as the
$\CC(\om)$-linear span of the $m$-diagrams with this
multiplication.

For $1\leqslant a<b\leqslant m$ denote by $s_{a\tss b}$
and $\ep_{a\tss b}$ the respective diagrams
of the form

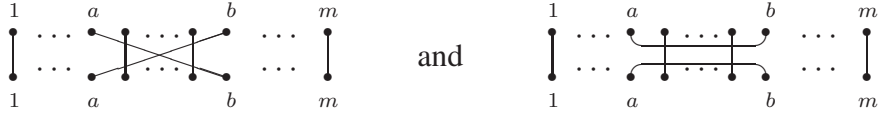
\setlength{\unitlength}{0.85pt}
\begin{center}
\begin{picture}(400,60)
\thinlines

\put(10,20){\circle*{3}}
\put(45,20){\circle*{3}}
\put(60,20){\circle*{3}}
\put(90,20){\circle*{3}}
\put(105,20){\circle*{3}}
\put(150,20){\circle*{3}}

\put(10,40){\circle*{3}}
\put(45,40){\circle*{3}}
\put(60,40){\circle*{3}}
\put(90,40){\circle*{3}}
\put(105,40){\circle*{3}}
\put(150,40){\circle*{3}}

\put(10,20){\line(0,1){20}}
\put(60,20){\line(0,1){20}}
\put(45,20){\line(3,1){60}}
\put(45,40){\line(3,-1){60}}
\put(90,20){\line(0,1){20}}
\put(150,20){\line(0,1){20}}

\put(20,20){$\cdots$}
\put(20,35){$\cdots$}
\put(68,20){$\cdots$}
\put(68,35){$\cdots$}
\put(120,20){$\cdots$}
\put(120,35){$\cdots$}

\put(8,47){\scriptsize $1$ }
\put(43,47){\scriptsize $a$ }
\put(105,47){\scriptsize $b$ }
\put(146,47){\scriptsize $m$ }

\put(8,7){\scriptsize $1$ }
\put(43,7){\scriptsize $a$ }
\put(105,7){\scriptsize $b$ }
\put(146,7){\scriptsize $m$ }

\put(190,25){\text{and}}

\put(250,20){\circle*{3}}
\put(285,20){\circle*{3}}
\put(300,20){\circle*{3}}
\put(330,20){\circle*{3}}
\put(345,20){\circle*{3}}
\put(390,20){\circle*{3}}

\put(250,40){\circle*{3}}
\put(285,40){\circle*{3}}
\put(300,40){\circle*{3}}
\put(330,40){\circle*{3}}
\put(345,40){\circle*{3}}
\put(390,40){\circle*{3}}

\put(250,20){\line(0,1){20}}
\put(300,20){\line(0,1){20}}
\put(315,20){\oval(60,12)[t]}
\put(315,40){\oval(60,12)[b]}
\put(330,20){\line(0,1){20}}
\put(390,20){\line(0,1){20}}

\put(260,20){$\cdots$}
\put(260,35){$\cdots$}
\put(308,20){$\cdots$}
\put(308,35){$\cdots$}
\put(360,20){$\cdots$}
\put(360,35){$\cdots$}

\put(248,47){\scriptsize $1$ }
\put(283,47){\scriptsize $a$ }
\put(345,47){\scriptsize $b$ }
\put(386,47){\scriptsize $m$ }

\put(248,7){\scriptsize $1$ }
\put(283,7){\scriptsize $a$ }
\put(345,7){\scriptsize $b$ }
\put(386,7){\scriptsize $m$ }

\end{picture}
\end{center}
\setlength{\unitlength}{1pt}

\noindent
They generate the algebra $\Bc_m(\om)$.
Its subalgebra spanned over $\CC$
by the diagrams without horizontal edges will be identified
with the group algebra of the symmetric group
$\CC[\Sym_m]$ so that $s_{a\tss b}$
is identified with the transposition $(a\ts b)$.

We will use a special element $s^{(m)}\in\Bc_m(\om)$,
known as the {\em symmetrizer}. Several explicit expressions for $s^{(m)}$
are collected in \cite[Ch.~1]{m:so}; we will recall one of them, as appeared in \cite{hx:ts},
\beql{symdia}
s^{(m)}=\frac{1}{m!}\ts\sum_{r=0}^{\lfloor m/2\rfloor}
(-1)^r\binom{\om/2+m-2}{r}^{-1}
\sum_{d\in\Dc^{(r)}}d,
\eeq
where $\Dc^{(r)}\subset\Bc_m(\om)$ denotes the set of diagrams
which have exactly $r$ horizontal edges in the top (and hence in the bottom) row.
Since $\Dc^{(0)}=\Sym_m$, the element
\beql{hm}
h^{(m)}=\frac{1}{m!}\ts \sum_{d\in\Dc^{(0)}}d
\eeq
coincides with the symmetrizer in $\CC[\Sym_m]$ in \eqref{ha}.

For every $a\in\{1,\dots,m\}$ introduce the transposition
$\tra_a$ as the linear map
\ben
\tra_a:\Bc_m(\om)\to \Bc_m(\om),\qquad d\mapsto d^{\ts\tra_a},
\een
where the diagram $d^{\ts\tra_a}$ is obtained
from $d$ by swapping the $a$-th dots in the top and bottom rows, as the ends of edges.
In particular, $s_{a\tss b}^{\ts\tra_a}=\ep_{a\tss b}$
and $\ep_{a\tss b}^{\ts\tra_a}=s_{a\tss b}$.
Denote by $\Jr_m$ the subspace of $\Bc_m(\om)$ spanned by
all sums $d+d^{\ts\tra_a}$ with $d\in \Bc_m(\om)$ and $a=1,\dots,m$.
Note that if $\tau=\tra_{a_1}\circ\dots\circ\tra_{a_s}$ is the composition
of an odd number of distinct transpositions, then the sum $d+d^{\ts\tau}$
belongs to $\Jr_m$.
Introduce a rational function in $\om$ by
\ben
\ga_m(\om)=\frac{\om+m-2}{\om+2\tss m-2}.
\een

\ble\label{lem:brarel}
In the algebra $\Bc_m(\om)$ we have
\ben
\ga_{m}(\om)\tss s^{(m)}\equiv h^{(m)}\mod \Jr_{m}.
\een
\ele

\bpf
If $m$ is odd, then the claim is trivial because both $s^{(m)}$ and $h^{(m)}$ belong
to the subspace $\Jr_{m}$. Now suppose $m=2k$ is even.
We will start with the formula \eqref{symdia} for $s^{(2k)}$
and use an inductive procedure to apply a sequence of reductions modulo $\Jr_{2k}$
to eliminate all diagrams containing horizontal edges from the sum.
Any diagram $d$ containing an edge of the form $(a,a)$ belongs to $\Jr_{2k}$,
so that such diagrams can be ignored in the procedure.
As a first step, for each $r=0,1,\dots,k$ split the set of diagrams $\Dc^{(r)}$
into three subsets,
\beql{dsp}
\Dc^{(r)}=\Dc^{(r,-)}\cup \Dc^{(r,0)}\cup \Dc^{(r,+)},
\eeq
where $d\in\Dc^{(r,-)}$ if and only if the vertices $1$ in the top and bottom rows
are the ends of horizontal edges; $d\in\Dc^{(r,+)}$ if and only if
the vertices $1$ are the ends of different non-horizontal edges, and the remaining
diagrams belong to $\Dc^{(r,0)}$. In particular,
$\Dc^{(k)}=\Dc^{(k,-)}$.
It is clear by the application of the transposition $\tra_1$ that for $r\geqslant 0$
\ben
\sum_{d\in \Dc^{(r,0)}}\tss d\equiv 0\mod \Jr_{2k}\Fand
\sum_{d\in \Dc^{(r+1,-)}}\tss d\equiv -\sum_{d\in \Dc^{(r,+)}}\tss d\mod \Jr_{2k}.
\een
Taking into account the relation
\ben
\binom{\om/2+2k-2}{r}^{-1}+\binom{\om/2+2k-2}{r+1}^{-1}=
\frac{\om+4\tss k-2}{\om+4\tss k-4}\ts \binom{\om/2+2k-3}{r}^{-1},
\een
we can conclude from \eqref{symdia} that the reduction modulo $\Jr_{2k}$ yields the equivalence
\beql{secre}
\ga_{2k}(\om)\tss s^{(2k)}\equiv
\frac{\ga_{2k-2}(\om+2)}{(2k)!}\ts\sum_{r=0}^{k-1}
(-1)^r\binom{\om/2+2k-3}{r}^{-1}
\sum_{d\in\Dc^{(r,+)}}d.
\eeq
Note that the inverse binomial coefficients in this expression coincide
with those in \eqref{symdia} for $m=2k-2$ with the parameter $\om$ replaced with $\om+2$.

For the second step of the reduction, represent each set $\Dc^{(r,+)}$
as the union
\ben
\Dc^{(r,+)}=\bigcup_{a,b=2}^{k} \Dc^{(r,+)}_{a,b},
\een
where the subset $\Dc^{(r,+)}_{a,b}$ consists of the diagrams $d$
containing the (non-horizontal) edges $(1,a)$ with the dot $1$ in the top row, and $(1,b)$
with the dot $1$ in the bottom row.
Re-arrange expression \eqref{secre} to include an extra sum by writing
\ben
\sum_{d\in\Dc^{(r,+)}}d=\sum_{a,b=2}^k\ \sum_{d\in\Dc^{(r,+)}_{a,b}}d
\een
and changing the order of summation to take the external sum over $a$ and $b$.
If $a=b$, then we proceed by applying the same reduction modulo $\Jr_{2k}$
as in the first step, by ignoring the vertices $1$ and $a$ in the top and bottom rows.
If $a\ne b$, then
split the union of sets $\Dc^{(r,+)}_{a,b}\cup \Dc^{(r,+)}_{b,a}$ as in \eqref{dsp},
\ben
\Dc^{(r,+)}_{a,b}\cup \Dc^{(r,+)}_{b,a}
=\Dc^{(r,+,-)}_{\{a,b\}}\cup \Dc^{(r,+,0)}_{\{a,b\}}\cup \Dc^{(r,+,+)}_{\{a,b\}},
\een
where $d\in\Dc^{(r,+,-)}_{\{a,b\}}$ if and only if the remaining
vertices $a$ and $b$
are the ends of horizontal edges; $d\in\Dc^{(r,+,+)}_{\{a,b\}}$ if and only if
the remaining
vertices $a$ and $b$ are the ends of different non-horizontal edges, and the remaining
diagrams belong to $\Dc^{(r,+,0)}_{\{a,b\}}$. Similar to the first reduction step,
the application of the composition of transpositions
$\tra_1\circ\tra_a\circ\tra_b$ shows that for $r\geqslant 0$
\ben
\sum_{d\in \Dc^{(r,+,0)}_{\{a,b\}}}\tss d\equiv 0\mod \Jr_{2k}\Fand
\sum_{d\in \Dc^{(r+1,+,-)}_{\{a,b\}}}\tss d\equiv -\sum_{d\in \Dc^{(r,+,+)}_{\{a,b\}}}\tss d\mod \Jr_{2k}.
\een
This leads to the second step reduction formula analogous to \eqref{secre},
and the argument continues in the same way, where compositions of $2r-1$ distinct transpositions
are used at the $r$-th step. As a result of the $k$-th step of the reduction procedure,
we get the sum of diagrams without horizontal edges with the overall coefficient $1/(2k)!$,
as required.
\epf

The Brauer algebra $\Bc_m(\om)$ with the special values $\om=N$ and $\om=-2n$
acts on the tensor space \eqref{tenprk}
so that
the action centralizers the respective diagonal actions of
the orthogonal and symplectic groups.
In the orthogonal case, the generators of $\Bc_m(N)$ act by the rule
\beql{braact}
s_{a\tss b}\mapsto P_{a\tss b},\qquad \ep_{a\tss b}\mapsto Q_{a\tss b},
\qquad 1\leqslant a<b\leqslant m,
\eeq
where $P_{a\tss b}$ is defined by \eqref{pdef},
while
\ben
Q_{a\tss b}=\sum_{i,j=1}^N 1^{\ot(a-1)}\ot e_{ij}
\ot 1^{\ot(b-a-1)}\ot e_{i'j'}\ot 1^{\ot(m-b)}.
\een
In the symplectic case, the action of
$\Bc_m(-N)$ with $N=2n$ in the space \eqref{tenprk}
is defined by
\beql{braactsp}
s_{a\tss b}\mapsto -P_{a\tss b},\qquad \ep_{a\tss b}\mapsto -Q_{a\tss b},\qquad
1\leqslant a<b\leqslant m,
\eeq
where
\ben
Q_{a\tss b}=\sum_{i,j=1}^{2n} \ve_i\tss\ve_j\ts 1^{\ot(a-1)}\ot e_{ij}
\ot 1^{\ot(b-a-1)}\ot e_{i'j'}\ot 1^{\ot(m-b)}.
\een

We will denote by
$S^{(m)}$ the image of the symmetrizer
$s^{(m)}\in\Bc_m(\om)$ under the respective actions \eqref{braact}
and \eqref{braactsp}, assuming $m\leqslant n$ in the symplectic case.
The elements
$\phi^{}_{m\tss m}$ of the vacuum module
$V_{\cri}(\g)\cong\U\big(t^{-1}\g[t^{-1}]\big)$ are defined
by
\beql{defphi}
\phi^{}_{m\tss m}=\ga_m(\om)\ts\tr^{}_{1,\dots,m}\ts S^{(m)}
\big(T+F[-1]_1\big)\dots \big(T+F[-1]_m\big)\tss 1
\eeq
where
$\om=N$ and $\om=-2n$, respectively, for the orthogonal and symplectic case.
In the symplectic case
the values of $m$ are restricted to $1\leqslant m\leqslant 2n+1$
with an additional justification of formula \eqref{defphi}
for the values $n+1\leqslant m\leqslant 2n+1$ via an `analytic continuation' argument;
see \cite[Sec.~8.3]{m:so}.
As proved in \cite{m:ff} (see also \cite[Ch.~8]{m:so}),
all elements $\phi^{}_{m\tss m}$ belong to the Feigin--Frenkel center $\z(\wh\g)$.
Moreover,
the elements $\phi^{}_{22},\phi^{}_{44},\dots,\phi^{}_{2n\tss 2n}$ form a complete
set of Segal--Sugawara vectors for $\g=\oa_{2n+1}$ and $\spa_{2n}$, whereas
$\phi^{}_{22},\phi^{}_{44},\dots,\phi^{}_{2n-2\ts 2n-2}, \Pf\tss F[-1]$ form a complete
set of Segal--Sugawara vectors for $\g=\oa_{2n}$.

Now we proceed in the same way as in the proof of Theorem~\ref{thm:typea}
by expanding the product in \eqref{defphi}
with the use of \cite[Lemmas~8.1.5 \& 8.3.1]{m:so} to get
\beql{expbcd}
\phi^{}_{m\tss m}=\ga_m(\om)\ts\tr^{}_{1,\dots,m}\ts S^{(m)} \sum_{\la\ts\vdash m}\tss
c^{}_{\la}\ts F[-\la],
\eeq
where we set
\ben
F[-\la]=F[-\la_1]_1 \dots F[-\la_{\ell}]_{\ell}
\een
for $\la=(\la_1,\dots,\la_{\ell})$ with $\ell=\ell(\la)$.
Note that the summands in \eqref{expbcd} with odd values of $\ell$ are equal to zero
because the matrices $F[-\la_i]$ are skew-symmetric
with respect to the transposition~$\tra$,
while $S^{(m)}$ is stable under the simultaneous transpositions
with respect to all $m$ copies of $\End\CC^N$.
Now use \cite[Eqs~(5.32) \& (5.43)]{m:so} to calculate partial traces to get
\beql{partsym}
\ga_m(-2n)\ts\tr^{}_{\ell+1,\dots,m}\ts S^{(m)}
=\binom{2n+1}{m}\binom{2n+1}{\ell}^{-1}\ga_{\ell}(-2n)\ts S^{(\ell)}
\eeq
in the symplectic case, and
\beql{partsymo}
\ga_m(N)\ts\tr^{}_{\ell+1,\dots,m}\ts S^{(m)}
=\binom{N+m-2}{m}\binom{N+\ell-2}{\ell}^{-1}\ga_{\ell}(N)\ts S^{(\ell)}
\eeq
in the orthogonal case.
The desired formulas \eqref{phipsibcd} are now implied by the relation
\beql{gatr}
\ga_{\ell}(\om)\ts\tr^{}_{1,\dots,\ell}\ts S^{(\ell)} F[-\la]
=\tr^{}_{1,\dots,\ell}\ts H^{(\ell)} F[-\la],
\eeq
where
$H^{(\ell)}$ denotes
the image of the element $h^{(\ell)}$ defined in \eqref{hm},
under the respective actions \eqref{braact}
and \eqref{braactsp} of the Brauer algebra.
Relation \eqref{gatr} follows from Lemma~\ref{lem:brarel}
because the transposition $\tra_a$ on the Brauer algebra is
consistent with the partial transposition $\tra_a$
on the tensor product \eqref{tenprkao}. That is,
if an element $s\in \Bc_{m}(\om)$ has the form $s=d+d^{\ts\tra_a}$, then
for the image $S$ of $s$ with respect to the actions \eqref{braact}
and \eqref{braactsp} we have
$\tr^{}_{1,\dots,m}\ts S\ts F[-\la]=0$. Indeed, one verifies easily
that $S$ is stable under the transposition $\tra_a$ on the space
\eqref{tenprkao}, so that the relation
holds since
$F[-\la_a]+F[-\la_a]^{\tra}=0$. It remains to note that
\beql{dpla}
D(\la)=\tr^{}_{1,\dots,\ell}\ts H^{(\ell)} F[-\la]\Fand
P(\la)=\tr^{}_{1,\dots,\ell}\ts H^{(\ell)} F[-\la]
\eeq
in the symplectic and orthogonal case, respectively.
\epf

\bre\label{rem:genco}
All coefficients of the polynomial in $T$ given by
\ben
\ga_m(\om)\ts\tr^{}_{1,\dots,m}\ts S^{(m)}
\big(T+F[-1]_1\big)\dots \big(T+F[-1]_m\big)=
\phi_{m\tss  0}\ts T^{\tss m}+\phi_{m\tss  1}\ts T^{\tss m-1}+\dots+\phi_{m\tss m}
\een
belong to $\z(\wh\g)$; see \cite[Ch.~8]{m:so}. By replacing $T\mapsto T+u$
for a variable $u$, one can show that
\ben
\phi_{m\tss  k}=\binom{N+m-2}{m-k}\tss \phi_{k\tss  k}\Fand
\phi_{m\tss  k}=\binom{2n-k+1}{m-k}\tss \phi_{k\tss  k}
\een
in the orthogonal and symplectic case, respectively.
\qed
\ere

\subsection{Symmetrization map}
\label{subsec:sm}

In her recent work \cite{o:sf}, Yakimova gave new formulas
for Segal--Sugawara vectors in types $B,C,D$ and $G_2$ by using
the canonical symmetrization map.
We will show that these vectors in the classical types coincide with
those found in \cite{m:ff}.

Recall that for a Lie algebra $\agot$ the symmetrization map
$\vpi:\Sr(\agot)\to\U(\agot)$ is defined by
\beql{sym}
\vpi:x_1\dots x_n\mapsto \frac{1}{n!}\sum_{\si\in\Sym_n}x_{\si(1)}\dots x_{\si(n)},\qquad x_i\in\agot.
\eeq
Regarding $T=-d/dt$ as a derivation of the Lie algebra $t^{-1}\g[t^{-1}]$, we will
apply the map $\vpi$ for the Lie algebra $\CC T\oplus t^{-1}\g[t^{-1}]$.

For any element $S\in \Sr(\g)$ we will denote by $S[-1]$ the image
of $S$ under the embedding $\Sr(\g)\hra\Sr\big(t^{-1}\g[t^{-1}]\big)$
defined by $X\mapsto X[-1]$ for $X\in\g$.

\paragraph{Type $A$.}
Take $\g=\gl_N$ and introduce elements $\De_k$ and $\Phi_k$ of the symmetric algebra $\Sr(\gl_N)$
by the expansions
\ben
\det(u+E)=u^N+\De_1\tss u^{N-1}+\dots+\De_N
\een
and
\ben
\det(1-q\tss E)^{-1}=1+\sum_{k=1}^{\infty} \Phi_k\ts q^k,
\een
for the matrix $E=[E_{ij}]$. Then
\beql{symingen}
\Sr(\gl_N)^{\gl_N}=\CC[\De_1,\dots,\De_N]=\CC[\Phi_1,\dots,\Phi_N].
\eeq

The first formula in the next proposition (along with its closely related
versions) was pointed out in \cite{o:sf}.

\bpr\label{prop:symal}
The Segal--Sugawara vectors \eqref{deftra} and \eqref{deftrh}
can be written in the form
\begin{align}\label{phimm}
\phi^{}_{m\tss m}&=\sum_{k=1}^m \binom{N-k}{m-k}\ts
\varpi\big(T^{m-k}\De_{k}[-1]\big)\tss 1\\
\intertext{and}
\label{psimm}
\psi^{}_{m\tss m}&=\sum_{k=1}^m \binom{N+m-1}{m-k}\ts
\varpi\big(T^{m-k}\tss\Phi_{k}[-1]\big)\tss 1.
\end{align}
\epr

\bpf
Expand the product in \eqref{deftra} as
\ben
\tr^{}_{1,\dots,m}\ts A^{(m)} \sum_{k=1}^m\
\sum_{1\leqslant i_1<\dots<i_{m-k}\leqslant m}E[-1]_1\dots
E[-1]_{i_1-1}T\tss E[-1]_{i_1+1}\dots E[-1]_m\tss 1,
\een
so that the factors $T$ occur in the places $i_1,\dots,i_{m-k}$. Now apply conjugations
by suitable elements of $\Sym_m$ and use the cyclic property of trace
to bring this expression to the form where the labels of the factors $E[-1]$
take consecutive values:
\begin{multline}
\non
\tr^{}_{1,\dots,m}\ts A^{(m)} \sum_{k=1}^m\
\sum_{1\leqslant j_1\leqslant \dots\leqslant j_{m-k}\leqslant k}E[-1]_1\dots
E[-1]_{j_1-1}T\tss E[-1]_{j_1}\\
{}\times\ldots\times E[-1]_{j_{m-k}-1}T\tss E[-1]_{j_{m-k}}\dots E[-1]_k\tss 1
\end{multline}
with the factors $T$ occurring just before $E[-1]_{j_a}$ for $a=1,\dots,m-k$.
Now apply formula \eqref{patr} for the partial traces of $A^{(m)}$
with $\ell$ replaced by $k$ to come to the expression
\ben
\sum_{k=1}^m\ts\frac{(N-k)!\ts k!}{(N-m)!\ts m!}\ts
\tr^{}_{1,\dots,k}\ts A^{(k)}
\sum_{1\leqslant j_1\leqslant \dots\leqslant j_{m-k}\leqslant k}E[-1]_1\dots
E[-1]_{j_1-1}T\tss E[-1]_{j_1}\dots E[-1]_k\tss 1.
\een
For a fixed value of $k$ we have the relation
\begin{multline}
\tr^{}_{1,\dots,k}\ts A^{(k)}
\sum_{1\leqslant j_1\leqslant \dots\leqslant j_{m-k}\leqslant k}E[-1]_1\dots
E[-1]_{j_1-1}T\tss E[-1]_{j_1}\dots E[-1]_k\tss 1\\
=\binom{m}{k}\ts\varpi\big(T^{m-k}\De_{k}[-1]\big)\tss 1.
\label{antsy}
\end{multline}
This proves \eqref{phimm}. The proof of \eqref{psimm} is the same, with the use
of \eqref{patrh}.
\epf

\paragraph{Type $C$.}
Write the elements $F_{ij}$ of the symplectic
Lie algebra $\spa_{2n}$ into the matrix $F=[F_{ij}]$.
Introduce elements $\De_{2l}$ of the symmetric algebra $\Sr(\spa_{2n})$
by
\ben
\det(u+F)=u^{2n}+\De_2\ts u^{2n-2}+\dots+\De_{2n}.
\een
We have
\ben
\Sr(\spa_{2n})^{\spa_{2n}}=\CC[\De_2,\De_4,\dots,\De_{2n}].
\een

\bpr\label{prop:symalc}
The Segal--Sugawara vectors \eqref{defphi}
can be written in the form
\beql{binc}
\phi^{}_{m\tss m}=\sum_{l=1}^{\lfloor m/2\rfloor} \binom{2n-2l+1}{m-2l}\ts
\varpi\big(T^{m-2l}\De_{2l}[-1]\big)\ts 1
\eeq
for $m=1,\dots,2n+1$.
\epr

\bpf
As in the proof of Theorem~\ref{thm:typebcd}, it will be sufficient to assume
that $m\leqslant n$. The arguments used in \cite[Sec.~8.3]{m:so}
will then allow one to extend the result to all remaining values of $m$.
Expand the product in \eqref{defphi} and apply conjugations by suitable
permutations to get
\beql{expbcdcomp}
\phi^{}_{m\tss m}=\ga_{m}(-2n)\ts\tr^{}_{1,\dots,m}\ts S^{(m)} \sum_{\be}\tss
d^{}_{\be}\ts F[-\be],
\eeq
summed over compositions $\be=(\be_1,\dots,\be_{\ell})$ of $m$ with $\be_i\geqslant 1$,
where we set
\ben
F[-\be]=F[-\be_1]_1 \dots F[-\be_{\ell}]_{\ell},
\een
while $d^{}_{\be}$ are certain integer coefficients.
As with the expansion \eqref{expbcd}, the summands with odd values of $\ell$ are equal to zero.
Now calculate partial traces by using \eqref{partsym} and apply relation
\eqref{gatr}, which holds in the same form for $F[-\la]$ replaced by $F[-\be]$, to get
\ben
\phi^{}_{m\ts m}=\sum_{\be}\tss
d^{}_{\be}\tss\binom{2n-2l+1}{m-2l}\binom{m}{2l}^{-1}\ts\tr^{}_{1,\dots,2l}\ts H^{(2l)}
F[-\be],
\een
summed over the compositions $\be=(\be_1,\dots,\be_{2l})$. As with relation \eqref{antsy},
for a fixed value of $l$ we have
\ben
\sum_{\be}\tss
d^{}_{\be}\ts\tr^{}_{1,\dots,2l}\ts H^{(2l)} F[-\be]=
\binom{m}{2l} \ts \varpi\big(T^{m-2l}\De_{2l}[-1]\big)\tss 1,
\een
thus competing the proof.
\epf

Proposition~\ref{prop:symalc} shows
that $\phi^{}_{m\tss m}$ with $m=2k$ coincides with the Segal--Sugawara vector
produced in \cite[Theorem~4.4]{o:sf}.

\paragraph{Types $B$ and $D$.}
Introduce elements $\Phi_{2l}$ of the symmetric algebra $\Sr(\oa_{N})$
by
\ben
\det(1-q\tss F)^{-1}=1+\sum_{k=1}^{\infty} \Phi_{2l}\ts q^{2l}
\een
for the matrix $F=[F_{ij}]$. Then
\ben
\Sr(\oa_{2n+1})^{\oa_{2n+1}}=\CC[\Phi_2,\dots,\Phi_{2n}]\Fand
\Sr(\oa_{2n})^{\oa_{2n}}=\CC\big[\Phi_2,\dots,\Phi_{2n-2},\Pf\tss F[-1]\big].
\een

\bpr\label{prop:symalbd}
The Segal--Sugawara vectors \eqref{defphi}
can be written in the form
\beql{bin}
\phi^{}_{m\tss m}=\sum_{l=1}^{\lfloor m/2\rfloor} \binom{N+m-2}{m-2l}
\varpi\big(T^{m-2l}\tss\Phi_{2l}[-1]\big)\ts 1
\eeq
for $m\geqslant 1$.
\epr

\bpf
The argument is the same as for Proposition~\ref{prop:symalc}, where
we use the partial trace formula \eqref{partsymo} instead of \eqref{partsym},
and the corresponding version of relation \eqref{gatr} for compositions.
\epf

Proposition~\ref{prop:symalbd} implies that
$\phi^{}_{m\tss m}$ with $m=2k$
coincides with the Segal--Sugawara vector
given by \cite[Theorem~7.6]{o:sf}, because the binomial coefficient
in \eqref{bin} coincides with the expression $R(k,k-l)$ used therein.

Comparing
the binomial coefficients in \eqref{binc} and \eqref{bin}
with those occurring in \eqref{patr} and \eqref{patrh},
one can write an equivalent
formula for the Segal--Sugawara vector $\phi^{}_{m\tss m}$ in the symbolic form
\beql{phisym}
\phi^{}_{m\tss m}=\tr^{}_{1,\dots,m}\ts \overline H^{\ts(m)}
\big(T+F[-1]_1\big)\dots \big(T+F[-1]_{m}\big)\tss 1,
\eeq
where the symbol $\overline H^{\ts(m)}$ is interpreted
as a modified operator $H^{(m)}$ via an expansion
similar to that of \eqref{defphi},
with the partial traces calculated by the rules
\ben
\tr^{}_m \overline H^{\ts(m)}=\frac{N+m-2}{m}\ts \overline H^{\ts(m-1)}\Fand
\tr^{}_m \overline H^{\ts(m)}=\frac{2n-m+2}{m}\ts \overline H^{\ts(m-1)}
\een
in the orthogonal and symplectic case, respectively.

To derive a formal expression in the symplectic case, introduce the $(2n+1)\times(2n+1)$
matrix $F^{\tss\circ}$ by inserting the zero row and column in the middle
of the matrix $F=[F_{ij}]$, as illustrated:
\beql{fcirc}
F^{\tss\circ}=\begin{bmatrix}
F_{11}&\dots& F_{1 n} & 0 & F_{1 n'} &\dots & F_{1 1'}\\
\vdots&\vdots&\vdots&\vdots&\vdots&\vdots&\vdots\\
F_{n1}&\dots& F_{n n} & 0 & F_{n n'} &\dots & F_{n 1'}\\
0&\cdots&0&0&0&\cdots&0\\
F_{n'1}&\dots& F_{n' n} & 0 & F_{n' n'} &\dots & F_{n' 1'}\\
\vdots&\vdots&\vdots&\vdots&\vdots&\vdots&\vdots\\
F_{1'1}&\dots& F_{1' n} & 0 & F_{1' n'} &\dots & F_{1' 1'}
\end{bmatrix}.
\eeq
Define also the corresponding matrices $F^{\tss\circ}[r]=\big[F^{\tss\circ}_{ij}[r]\big]$
and label their rows and columns by the symbols
$1,\dots,n,0,n'\dots,1'$ which will also label the canonical basis elements
of the space $\CC^{2n+1}$. Consider the tensor space \eqref{tenprkao} for $\g=\spa_{2n}$
with $N=2n+1$
and let $A^{(m)}$ be the
anti-symmetrization operator in this space, as defined in Sec.~\ref{subsec:ga}.
Now introduce elements $\phi_{m\tss a}\in V_{\cri}(\spa_{2n})\cong\U\big(t^{-1}\spa_{2n}[t^{-1}]\big)$
by the expansion
\beql{traci}
\tr^{}_{1,\dots,m}\ts A^{(m)}
\big(T+F^{\tss\circ}[-1]_1\big)\dots \big(T+F^{\tss\circ}[-1]_m\big)
=\phi_{m\tss  0}\ts T^{\tss m}+\phi_{m\tss  1}\ts T^{\tss m-1}+\dots+\phi_{m\tss m}.
\eeq

\bco\label{cor:nnpo}
All elements $\phi_{m\tss a}$ belong to the Feigin--Frenkel center
$\z(\wh\spa_{2n})$.
Moreover, the family
$\phi^{}_{2\tss 2},\phi^{}_{4\tss 4},\dots,\phi^{}_{2n\ts 2n}$
is a complete set of Segal--Sugawara vectors.
\eco

\bpf
Similar to \eqref{expbcdcomp}, expand the product on the left hand side of \eqref{traci}
to bring it to the form
\beql{afcir}
\tr^{}_{1,\dots,m}\ts A^{(m)} \sum_{k=0}^m\Big(\sum_{\be}\tss
d^{\tss(m)}_{\be}\ts F^{\tss\circ}[-\be]\Big)\ts T^{\tss m-k}
\eeq
with the second sum over compositions $\be=(\be_1,\dots,\be_{\ell})$ of $k$ with $\be_i\geqslant 1$,
where
\ben
F^{\tss\circ}[-\be]=F^{\tss\circ}[-\be_1]_1 \dots F^{\tss\circ}[-\be_{\ell}]_{\ell},
\een
and $d^{\tss(m)}_{\be}$ are certain integer coefficients.
Now use \eqref{patr} to write \eqref{afcir} as
\beql{expfi}
\sum_{k=0}^m\Big(\sum_{\be}\tss d^{\tss(m)}_{\be}\ts
\binom{2n+1}{m}\binom{2n+1}{\ell}^{-1} \tr^{}_{1,\dots,\ell}\ts A^{(\ell)}
F^{\tss\circ}[-\be]\Big)\ts T^{\tss m-k}.
\eeq
Now observe that $\tr^{}_{1,\dots,\ell}\ts A^{(\ell)}F^{\tss\circ}[-\be]$
coincides with the trace $\tr^{}_{1,\dots,\ell}\ts A^{(\ell)}F^{}[-\be]$
calculated over the tensor space \eqref{tenprkao} for $\g=\spa_{2n}$
with $N=2n$.

On the other hand, the arguments used in the proof of Theorem~\ref{thm:typebcd}
show that the polynomial in $T$ defined by
\beql{est}
\ga_m(-2n)\ts\tr^{}_{1,\dots,m}\ts S^{(m)}
\big(T+F[-1]_1\big)\dots \big(T+F[-1]_m\big)
\eeq
is given by the same expression \eqref{expfi}. Hence, the polynomials \eqref{traci}
and \eqref{est} coincide, and
the corollary follows from the results of \cite[Sec.~8.3]{m:so}.
\epf

\section{Casimir elements and Harish-Chandra images}
\label{sec:ce}

\subsection{Symmetrized basic invariants}
\label{subsec:sb}

\paragraph{Type $A$.}
Applying the symmetrization map
\eqref{sym} for the Lie algebra $\gl_N$, and using \eqref{symingen}, we get algebraically independent
generators of the center $\Zr(\gl_N)$ of the universal enveloping algebra $\U(\gl_N)$,
\ben
\Zr(\gl_N)=\CC\big[\vpi(\De_1),\dots,\vpi(\De_N)\big]=\CC\big[\vpi(\Phi_1),\dots,\vpi(\Phi_N)\big].
\een

Given an $N$-tuple of complex numbers
$\la=(\la_1,\dots,\la_N)$, the corresponding irreducible highest weight
representation
$L(\la)$ of $\gl_N$ is generated by
a nonzero vector $\xi\in L(\la)$
such that
\begin{alignat}{2}
E_{ij}\ts\xi&=0 \qquad &&\text{for} \quad
1\leqslant i<j\leqslant N, \qquad \text{and}
\non
\\
E_{ii}\ts\xi&=\la_i\ts\xi \qquad &&\text{for} \quad 1\leqslant i\leqslant N.
\non
\end{alignat}
Any element $z\in\Z(\gl_N)$ acts in $L(\la)$
by multiplying each vector by a scalar $\chi(z)$.
As a function of the highest weight, $\chi(z)$
is a {\em shifted symmetric polynomial} in the variables $\la_1,\dots,\la_N$
which can be regarded as the image of $z$ under the
Harish-Chandra isomorphism $\chi$. This polynomial
is symmetric in the shifted variables $\la_1,\la_2-1,\dots,\la_N-N+1$.

Consider the {\em elementary shifted symmetric polynomials}
\ben
e^*_m(\la_1,\dots,\la_N)=\sum_{i_1<\dots<i_m}\la_{i_1}(\la_{i_2}-1)\dots
(\la_{i_m}-m+1)
\een
and the
{\em complete shifted symmetric polynomials}
\ben
h^*_m(\la_1,\dots,\la_N)=\sum_{i_1\leqslant\dots\leqslant i_m}\la_{i_1}(\la_{i_2}+1)\dots
(\la_{i_m}+m-1).
\een
They are particular cases of the {\em factorial} or {\em shifted Schur polynomials}; see \cite{m:sft}
and \cite{oo:ss}.

Recall that the {\em Stirling number of the second kind} ${\displaystyle \stir{m}{k}}$ counts the number
of partitions
of the set $\{1,\dots,m\}$ into $k$ nonempty subsets.

\bth\label{thm:chim} For the Harish-Chandra images we have
\begin{align}
\chi:\vpi(\De_m)&\mapsto\sum_{k=1}^m\stir{m}{k}\binom{N}{m}\binom{N}{k}^{-1}\ts
e^*_k(\la_1,\dots,\la_N)
\non\\
\intertext{and}
\chi:\vpi(\Phi_m)&\mapsto\sum_{k=1}^m\stir{m}{k}\binom{-N}{m}\binom{-N}{k}^{-1}\ts
h^*_k(\la_1,\dots,\la_N).
\non
\end{align}
\eth

\bpf
We will use the matrix notation of Sec.~\ref{subsec:ga}
applied to the algebra $\U(\gl_N)$ in place of $\U(\wh\gl_N)$.
Regarding the matrix $E=[E_{ij}]$ as the element
\ben
E=\sum_{i,j=1}^N e_{ij}\ot E_{ij}\in \End\CC^{N}\ot \U(\gl^{}_N)
\een
we get the following counterparts of \eqref{deftra} and \eqref{deftrh}:
\ben
\vpi(\De_m)=\tr^{}_{1,\dots,m}\ts A^{(m)}E_1\dots E_m\Fand
\vpi(\Phi_m)=\tr^{}_{1,\dots,m}\ts H^{(m)}E_1\dots E_m.
\een
On the other hand, the Harish-Chandra images
\ben
\chi:\tr^{}_{1,\dots,m}\ts A^{(m)}E_1(E_2-1)\dots (E_m-m+1)\mapsto e^*_m(\la_1,\dots,\la_N)
\een
and
\ben
\chi:\tr^{}_{1,\dots,m}\ts H^{(m)}E_1(E_2+1)\dots (E_m+m-1)\mapsto h^*_m(\la_1,\dots,\la_N)
\een
are well-known; see e.g. \cite[Secs~4.6 \& 4.7]{m:so} for proofs.
By the same argument as used in the proof of Theorem~\ref{thm:typea}, the identity
\beql{stide}
x^m=\sum_{k=1}^m\stir{m}{k}\ts x(x-1)\dots (x-k+1)
\eeq
implies that
\ben
\tr^{}_{1,\dots,m}\ts A^{(m)}E_1\dots E_m=\tr^{}_{1,\dots,m}\ts A^{(m)}\sum_{k=1}^m\stir{m}{k}\ts
E_1(E_2-1)\dots (E_k-k+1)
\een
and
\ben
\tr^{}_{1,\dots,m}\ts H^{(m)}E_1\dots E_m=\tr^{}_{1,\dots,m}\ts H^{(m)}\sum_{k=1}^m(-1)^{m-k}\ts
\stir{m}{k}\ts
E_1(E_2+1)\dots (E_k+k-1).
\een
The required formulas now follow by calculating
the partial traces over the spaces $\End\CC^{N}$ labelled by $k+1,\dots,m$, with the
use of \eqref{patr} and \eqref{patrh}.
\epf

\paragraph{Types $B$, $C$ and $D$.}
Now use the notation of Secs~\ref{subsec:nf} \& \ref{subsec:sm} and let $\g$ be
the orthogonal Lie algebra $\oa_N$ with $N=2n+1$ and $N=2n$
or symplectic Lie algebra $\spa_N$ with $N=2n$.

Given any $n$-tuple of complex numbers
$\la=(\la_1,\dots,\la_n)$, the corresponding irreducible highest weight
representation
$L(\la)$ of the Lie algebra $\g$ is generated by
a nonzero vector $\xi\in L(\la)$
such that
\begin{alignat}{2}
F_{ij}\ts\xi&=0 \qquad &&\text{for} \quad
1\leqslant i<j\leqslant N, \qquad \text{and}
\non
\\
F_{ii}\ts\xi&=\la_i\ts\xi \qquad &&\text{for} \quad 1\leqslant i\leqslant n.
\non
\end{alignat}
Any element $z$ of the center of $\U(\g)$ acts in $L(\la)$
by multiplying each vector by a scalar $\chi(z)$.
As a function of the highest weight, $\chi(z)$
is a {\em shifted invariant polynomial} in the variables $\la_1,\dots,\la_n$
with respect to the action of the corresponding Weyl group. The polynomial $\chi(z)$
can be regarded
as the Harish-Chandra image of $z$.

\bth\label{thm:hchbcd}
\begin{enumerate}
\item
If $\g=\spa_{2n}$, then for $m=2,4,\dots,2n$ the Harish-Chandra images are
\ben
\chi:\vpi(\De_m)\mapsto\sum_{k=1}^m\stir{m}{k}\binom{2n+1}{m}\binom{2n+1}{k}^{-1}
\ts e^*_k(\la_1,\dots,\la_n,0,-\la_n,\dots,-\la_1).
\een
\item
If $\g=\oa_{2n+1}$, then for even $m\geqslant 2$ the Harish-Chandra images are
\ben
\chi:\vpi(\Phi_m)\mapsto\sum_{k=1}^m\stir{m}{k}\binom{-2n}{m}\binom{-2n}{k}^{-1}
\ts h^*_k(\la_1,\dots,\la_n,-\la_n,\dots,-\la_1).
\een
\item
If $\g=\oa_{2n}$, then for even $m\geqslant 2$ the Harish-Chandra images are
\begin{multline}
\chi:\vpi(\Phi_m){}\mapsto\sum_{k=1}^m\stir{m}{k}\binom{-2n+1}{m}\binom{-2n+1}{k}^{-1}
\non\\[0.5em]
{}\times\ts
\Big(\frac12\ts h^*_k(\la_1,\dots,\la_{n-1},-\la_n,\dots,-\la_1)
+\frac12\ts h^*_k(\la_1,\dots,\la_n,-\la_{n-1},\dots,-\la_1)\Big).
\non
\end{multline}
\end{enumerate}
\eth

\bpf
The trace
\ben
\tr^{}_{1,\dots,m}\ts H^{(m)}F_1\dots F_m
\een
coincides with $\vpi(\De_m)$ and $\vpi(\Phi_m)$, respectively, in the symplectic
and orthogonal case. By applying Lemma~\ref{lem:brarel} as in the proof of Theorem~\ref{thm:typebcd},
we get
\ben
\tr^{}_{1,\dots,m}\ts H^{(m)}F_1\dots F_m=
\ga_m(\om)\ts\tr^{}_{1,\dots,m}\ts S^{(m)}
F_1\dots F_m,
\een
where
$\om=N$ and $\om=-2n$, respectively, for the orthogonal and symplectic case.
Using \eqref{stide} again, and adjusting the arguments of
the proof of Theorem~\ref{thm:typebcd} to the case of Lie algebra $\g$, we derive the relations
\ben
\ga_m(-2n)\ts\tr^{}_{1,\dots,m}\ts S^{(m)}
F_1\dots F_m=\ga_m(-2n)\ts\tr^{}_{1,\dots,m}\ts S^{(m)}\sum_{k=1}^m\stir{m}{k}\ts
F_1(F_2-1)\dots (F_k-k+1)
\een
in the symplectic case, and
\begin{multline}
\non
\ga_m(N)\ts\tr^{}_{1,\dots,m}\ts S^{(m)}
F_1\dots F_m\\
{}=\ga_m(N)\ts\tr^{}_{1,\dots,m}\ts S^{(m)}\sum_{k=1}^m(-1)^{m-k}\ts
\stir{m}{k}\ts
F_1(F_2+1)\dots (F_k+k-1)
\end{multline}
in the orthogonal case. By the results of \cite[Sec.~6]{mm:yc} (see also \cite[Sec.~13.4]{m:so}),
we have the Harish-Chandra images
\ben
\chi: \ga_k(-2n)\ts\tr^{}_{1,\dots,k}\ts S^{(k)}
F_1(F_2-1)\dots (F_k-k+1)\mapsto
e^*_k(\la_1,\dots,\la_n,0,-\la_n,\dots,-\la_1)
\een
for $\g=\spa_{2n}$,
\ben
\chi: \ga_k(N)\ts\tr^{}_{1,\dots,k}\ts S^{(k)}
F_1(F_2+1)\dots (F_k+k-1)\mapsto
h^*_k(\la_1,\dots,\la_n,-\la_n,\dots,-\la_1)
\een
for $\g=\oa_{2n+1}$, and
\begin{multline}
\chi: \ga_k(N)\ts\tr^{}_{1,\dots,k}\ts S^{(k)}
F_1(F_2+1)\dots (F_k+k-1)\\[0.4em]
{}\mapsto
\Big(\frac12\ts h^*_k(\la_1,\dots,\la_{n-1},-\la_n,\dots,-\la_1)
+\frac12\ts h^*_k(\la_1,\dots,\la_n,-\la_{n-1},\dots,-\la_1)\Big)
\non
\end{multline}
for $\g=\oa_{2n}$. The proof is completed by calculating
the partial traces of $\ga_m(\om)\ts S^{(m)}$
over the spaces $\End\CC^{N}$ labelled by $k+1,\dots,m$, with the
use of \eqref{partsym} and \eqref{partsymo}.
\epf

\bre\label{rem:mze}
As the proof of Theorem~\ref{thm:hchbcd} shows, the formulas for the Harish-Chandra
images extend to odd values of $m$, assuming that $\De_m=\Phi_m=0$. This provides
linear dependence relations for the elementary and complete shifted symmetric polynomials.
\qed
\ere

\subsection{Capelli-type determinant for $\spa_{2n}$}
\label{subsec:ct}

The coefficients of the polynomial in $u$ defined by
\begin{multline}
\non
D_{m}(u)=
\ga_m(\om)\ts\tr^{}_{1,\dots,m}\ts S^{(m)}\Big(F_1+u+\frac{m-1}{2}\Big)\Big(F_2+u+\frac{m-3}{2}\Big)\\
{}\times \dots \times\Big(F_{m}+u-\frac{m-1}{2}\Big),
\end{multline}
where $\om=N$ and $\om=-2n$ in the cases $\g=\oa_N$ and $\g=\spa_{2n}$, respectively,
are Casimir elements for $\g$. Their images under the Harish-Chandra isomorphism
were found in \cite{imr:ce}; see also \cite[Ch.~5]{m:so}. As with the Segal--Sugawara vectors,
the arguments of Sec.~\ref{subsec:sm} imply the interpretations for $D_{m}(u)$
in a symbolic form,
analogous to \eqref{phisym}:
\ben
D_{m}(u)=\tr^{}_{1,\dots,m}\ts \overline H^{\tss(m)}
\Big(F_1+u+\frac{m-1}{2}\Big)\Big(F_2+u+\frac{m-3}{2}\Big)
\dots \Big(F_{m}+u-\frac{m-1}{2}\Big).
\een

In the rest of this section we assume that $\g=\spa_{2n}$. Using the matrix \eqref{fcirc}
introduce the {\em symmetrized determinant} $C(u)$
of the matrix $u+F^{\tss\circ}:=u\tss I+F^{\tss\circ}$ by
\begin{multline}
\non
u\ts C(u)=\frac{1}{(2n+1)!}\ts\sum_{\si,\tau}\sgn\si\tau\cdot
\big(u+n+F^{\tss\circ}\big)_{\si(1)\tss\tau(1)}\big(u+n-1+F^{\tss\circ}\big)_{\si(2)\tss\tau(2)}\\
{}\times \dots \times\big(u-n+F^{\tss\circ}\big)_{\si(2n)\ts\tau(2n)},
\end{multline}
summed over all permutations $\si,\tau$ of the set $\{1,\dots,n,0,n',\dots,1'\}$.

\bpr\label{prop:capsp}
The symmetrized determinant is an even polynomial in $u$
of the form
\ben
C(u)=u^{2n}+C_1\ts u^{2n-2}+\dots+C_n,
\een
and the coefficients $C_1,\dots,C_n$ are algebraically independent generators
of the center of the algebra $\U(\spa_{2n})$.
Moreover, for the Harish-Chandra image we have
\ben
\chi:C(u)\mapsto (u^2-l_1^2)\dots (u^2-l_n^2),
\een
where $l_i=\la_i+n-i+1$ for $i=1,\dots,n$.
\epr

\bpf
Consider the tensor product algebra
\eqref{tenprkao} for $N=2n+1$ with the last tensor factor replaced by $\U(\spa_{2n})$.
As with Corollary~\ref{cor:nnpo},
we can write an equivalent expression for the polynomial $D_{m}(u)$ in the form
\ben
D_{m}(u)=\tr^{}_{1,\dots,m}\ts A^{\tss(m)}
\Big(F^{\tss\circ}_1+u+\frac{m-1}{2}\Big)\Big(F^{\tss\circ}_2+u+\frac{m-3}{2}\Big)
\dots \Big(F^{\tss\circ}_{m}+u-\frac{m-1}{2}\Big).
\een
Then clearly $u\ts C(u)=D^{}_{2n+1}(u)$ so that all claims follow from
\cite[Cor.~5.5.5]{m:so} and the identity
\ben
(u^2-l_1^2)\dots (u^2-l_n^2)=\sum_{k=0}^n(-1)^k\ts e_k(l_1^2,\dots, l_n^2\ts |\ts a)\ts
(u^2-1^2)\dots (u^2-(n-k)^2)
\een
for the factorial elementary symmetric polynomials $e_k(l_1^2,\dots, l_n^2\ts |\ts a)$
associated with the sequence $a=(1^2,2^2,\dots)$; see \cite[Eq.~(6.5)]{m:sft}.
\epf

Proposition~\ref{prop:capsp} provides a new formula for the Capelli-type determinant
introduced in \cite{m:sd}, which also plays the role of a noncommutative characteristic
polynomial of the matrix $F$.

%\newpage
\bigskip\bigskip

\small

\noindent
School of Mathematics and Statistics\newline
University of Sydney,
NSW 2006, Australia\newline
alexander.molev@sydney.edu.au

\end{document}